\documentclass[11pt]{article}
\usepackage{amsmath,amssymb}

\newtheorem{propo}{{\bf Proposition}}[section]
\newtheorem{coro}[propo]{{\bf Corollary}}
\newtheorem{lemma}[propo]{{\bf Lemma}} \newtheorem{theor}[propo]{{\bf
Theorem}} \newtheorem{ex}{{\sc Example}}[section]

\begin{document}

\vspace*{1.0in}

\begin{center} SOLVABLE LIE $A$-ALGEBRAS 
\end{center}
\bigskip

\begin{center} DAVID A. TOWERS 
\end{center}
\bigskip
\centerline {Department of
Mathematics, Lancaster University} \centerline {Lancaster LA1 4YF,
England} \centerline {Email: d.towers@lancaster.ac.uk}
\bigskip

\begin{abstract}
A finite-dimensional Lie algebra $L$ over a field $F$ is called an $A$-algebra if all of its nilpotent subalgebras are abelian. This is analogous to the concept of an $A$-group: a finite group with the property that all of its Sylow subgroups are abelian. These groups were first studied in the 1940s by Philip Hall, and are still studied today. Rather less is known about $A$-algebras, though they have been studied and used by a number of authors. The purpose of this paper is to obtain more detailed results on the structure of solvable Lie $A$-algebras.
\par
It is shown that they split over each term in their derived series. This leads to a decomposition of $L$ as $L = A_{n} \dot{+} A_{n-1} \dot{+} \ldots \dot{+} A_0$ where $A_i$ is an abelian subalgebra of $L$ and $L^{(i)} = A_{n} \dot{+} A_{n-1} \dot{+} \ldots \dot{+} A_{i}$ for each $0 \leq i \leq n$. It is shown that the ideals of $L$ relate nicely to this decomposition: if $K$ is an ideal of $L$ then $K = (K \cap A_n) \dot{+} (K \cap A_{n-1}) \dot{+} \ldots \dot{+} (K \cap A_0)$. When $L^2$ is nilpotent we can locate the position of the maximal nilpotent subalgebras: if $U$ is a maximal nilpotent subalgebra of $L$ then $U = (U \cap L^2) \oplus (U \cap C)$ where $C$ is a Cartan subalgebra of $L$.
\par 
If $L$ has a unique minimal ideal $W$ then $N = Z_L(W)$. If, in addition, $L$ is strongly solvable the maximal nilpotent subalgebras of $L$ are $L^2$ and the Cartan subalgebras of $L$ (that is, the subalgebras that are complementary to $L^2$.) Necessary and sufficient conditions are given for such an algebra to be an $A$-algebra. Finally, more detailed structure results are given when the underlying field is algebraically closed.  
\end{abstract}

\section{Introduction}
\medskip
A finite-dimensional Lie algebra $L$ over a field $F$ is called an $A$-algebra if all of its nilpotent subalgebras are abelian. This is analogous to the concept of an $A$-group, which is a finite group with the property that all of its Sylow subgroups are abelian. These groups were first studied in the 1940s by Philip Hall as soluble $A$-groups, and are still studied today. A great deal is known about their structure. Rather less is known about $A$-algebras, though they have been studied and used by a number of authors, including Bakhturin and Semenov \cite{bs}, Dallmer \cite{dall}, Drensky \cite{dren}, Sheina \cite{sh1} and \cite{sh2}, Premet and Semenov \cite{ps}, Semenov \cite{sem} and Towers and Varea \cite{tv1}, \cite{tv2}. 
\par
They arise in the study of constant Yang-Mills potentials. Every non-abelian nilpotent Lie algebra admits a non-trivial solution of the constant Yang-Mills equations. Moreover, if a subalgebra of a Lie algebra $L$ admits a non-trivial solution of the Yang-Mills equations then so does $L$. It is therefore useful to know if a given non-nilpotent Lie algebra has a non-abelian nilpotent subalgebra (see \cite{dall} for more details). They have also been particularly important in relation to the problem of describing residually finite varieties (see \cite{bs}, \cite{sh1}, \cite{sh2}, \cite{sem} and \cite{ps}).
\par
The Frattini ideal of $L$, $\phi(L)$, is the largest ideal of $L$ contained in
all maximal subalgebras of $L$. The Lie algebra $L$ is called {\it $\phi$-free}
if $\phi(L) = 0$, and {\it elementary} if $\phi(B)=0$ for every subalgebra $B$ of $L$. We say that $L$ is an $E$-algebra if $\phi(B)\leq \phi(L)$ for all subalgebras $B$ of $L$. Following Jacobson \cite{jac}, we say that a linear Lie algebra $L\leq {\rm gl}(V)$ is {\it almost algebraic} if $L$ contains the
nilpotent and semisimple Jordan components of its elements. Every algebraic Lie algebra is almost algebraic. An abstract Lie algebra $L$ is called almost algebraic if ${\rm ad}L\leq {\rm gl}(L)$ is almost algebraic. The classes of elementary Lie algebras, $E$-algebras, almost algebraic Lie algebras and $A$-algebras are related, as is shown in \cite{tv1} and \cite{tv2}. The {\em centre} of $L$ is $Z(L) = \{x \in L: [x,y] = 0 \hspace{.1cm} {\rm for \hspace{.1cm} all} \hspace{.1cm} y \in L \}$. We summarise below some of the known results for Lie $A$-algebras.
\bigskip

\begin{theor}\label{t:a}
Let $L$ be a Lie $A$-algebra over a field $F$.
\begin{itemize}
\item[(i)] If $F$ has characteristic zero, then
\begin{itemize}
\item[(a)] $L$ is almost algebraic if and only if it is elementary; in this case $L$ splits over each of its ideals;
\item[(b)] $L$ is elementary whenever L/R(L) and R(L) are elementary, where $R(L)$ is the solvable radical of $L$; and
\item[(c)] $L$ is an $E$-algebra.
\end{itemize}
\item[(ii)] If $F$ has characteristic $\neq 2, 3$, then $Q(L) = \{c \in L : {\rm (ad}c)^2 = 0 \}$is the unique maximal abelian ideal in $L$.
\item[(iii)] If $F$ has characteristic $\neq 2, 3$ and cohomological dimension $\leq 1$, then
\begin{itemize}
\item[(a)] $L^2 \cap Z(L) = 0$; and
\item[(b)] $L$ has a Levi decomposition and every Levi subalgebra is representable as a direct sum of simple ideals, each one of which splits over some finite extension of the ground field into a direct sum of ideals isomorphic to $sl(2)$.
\end{itemize}
\end{itemize}
\end{theor}
\medskip
{\it Proof.} (i) See Towers and Varea, \cite{tv2}.
\par
(ii), (iii) See Premet and Semenov, \cite{ps}.
\bigskip

The purpose of this paper is to obtain more detailed results on the structure of solvable Lie $A$-algebras. Some of the development is suggested by \cite{mak}, but more is possible for Lie algebras.
\par
In section two we collect together the preliminary results that we need, including the fact that for Lie $A$-algebras the derived series coincides with the lower nilpotent series. We also see that Lie $A$-algebras need not be metabelian.
\par
Section three contains the basic structure theorems for solvable Lie $A$-algebras. First they split over each term in their derived series. This leads to a decomposition of $L$ as $L = A_{n} \dot{+} A_{n-1} \dot{+} \ldots \dot{+} A_0$ where $A_i$ is an abelian subalgebra of $L$ and $L^{(i)} = A_{n} \dot{+} A_{n-1} \dot{+} \ldots \dot{+} A_{i}$ for each $0 \leq i \leq n$. It is shown that the ideals of $L$ relate nicely to this decomposition: if $K$ is an ideal of $L$ then $K = (K \cap A_n) \dot{+} (K \cap A_{n-1}) \dot{+} \ldots \dot{+} (K \cap A_0)$; moreover, if $N$ is the nilradical of $L$, $Z(L^{(i)}) = N \cap A_i$. We also see that the result in Theorem \ref{t:a} (iii)(a) holds when $L$ is solvable without any restrictions on the underlying field. 
\par
The fourth section looks at Lie $A$-algebras in which $L^2$ is nilpotent. These are metabelian and so the results of section three simplify. In addition we can locate the position of the maximal nilpotent subalgebras: if $U$ is a maximal nilpotent subalgebra of $L$ then $U = (U \cap L^2) \oplus (U \cap C)$ where $C$ is a Cartan subalgebra of $L$. 
\par
Section five is devoted to Lie $A$-algebras having a unique minimal ideal $W$. These have played a significant part in the study of varieties of residually finite Lie algebras. Again some of the results of sections three and four simplify. In particular, $N = Z_L(W)$, and if $L$ is strongly solvable the maximal nilpotent subalgebras of $L$ are $L^2$ and the Cartan subalgebras of $L$ (that is, the subalgebras that are complementary to $L^2$.) We also give necessary and sufficient conditions for a Lie algebra with a unique minimal ideal to be a strongly solvable $A$-algebra.  
\par
The final section is devoted to more detailed structure results when the underlying field is algebraically closed.
\par
Throughout $L$ will denote a finite-dimensional Lie algebra over a field $F$. Algebra direct sums will be denoted by $\oplus$, whereas vector space direct sums will be denoted by $\dot{+}$.
\bigskip

\section{Preliminary results}
\medskip

First we note that the class of Lie $A$-algebras is closed with respect to subalgebras, factor algebras and direct sums. Also that there is always a unique maximal abelian ideal, and it is the nilradical (which is equal to $Q(L)$ if $F$ has characteristic $\neq 2,3$, by Theorem \ref{t:a} (ii)).
\bigskip
 
\begin{lemma}\label{l:lemm1}
Let $L$ be a Lie $A$-algebra and let $N$ be its nilradical. Then
\begin{itemize}
\item[(i)] $N$ is the unique maximal abelian ideal of $L$;
\item[(ii)] if $B$ and $C$ are abelian ideals of $L$, we have $[B, C] = 0$; and
\item[(iii)] every subalgebra and every factor algebra of $L$ is an $A$-algebra.
\end{itemize}
\end{lemma}
\medskip
{\it Proof.} (i) Clearly $N$ is abelian and contains every abelian ideal of $L$.
\par

(ii) Simply note that $B, C \subseteq N$.
\par

(iii) It is easy to see that $L$ is subalgebra closed; that it is factor algebra closed is \cite[Lemma 1]{ps}.
\bigskip

\begin{lemma}\label{l:lemm2} Let $B$, $C$ be ideals of the Lie algebra $L$.
\begin{itemize}
\item[(i)] If $L/B$, $L/C$ are $A$-algebras, then $L/(B \cap C)$ is an $A$-algebra.
\item[(ii)] If $L = B \oplus C$, where $B, C$ are $A$-algebras, then $L$ is an $A$-algebra.
\end{itemize}
\end{lemma}
\medskip
{\it Proof.} (i) Let $U/(B \cap C)$ be a nilpotent subalgebra of $L/(B \cap C)$. Then $(U + B)/B$ is a nilpotent subalgebra of $L/B$, which is an $A$-algebra. It follows that $U^2 \subseteq B$. Similarly, $U^2 \subseteq C$, whence the result.
\par
(ii) This follows from (i).
\bigskip

We define the {\em nilpotent residual}, $\gamma_{\infty}(L)$, of $L$ be the smallest ideal of $L$ such that $L/\gamma_{\infty}(L)$ is nilpotent. Clearly this is the intersection of the terms of the lower central series for $L$. Then the {\em lower nilpotent series} for $L$ is the sequence of ideals $N_i(L)$ of $L$ defined by $N_0(L) = L$, $N_{i+1}(L) = \gamma_{\infty}(N_i(L))$ for $i \geq 0$. The {\em derived series} for $L$ is the sequence of ideals $L^{(i)}$ of $L$ defined by $L^{(0)} = L$, $L^{(i+1)} = [L^{(i)},L^{(i)}]$ for $i \geq 0$; we will also write $L^2$ for $L^{(1)}$. If $L^{(n)} = 0$ but $L^{(n-1)} \neq 0$ we say that that $L$ has {\em derived length} $n$.
\par
For Lie $A$-algebras we have the following result.
\bigskip

\begin{lemma}\label{l:series}
Let $L$ be a Lie $A$-algebra. Then the lower nilpotent series coincides with the derived series.
\end{lemma}
\medskip
{\it Proof.} Since $L/L^{(1)}$ is nilpotent we have $N_1(L) \subseteq L^{(1)}$. Also $L/N_1(L)$ is nilpotent and hence abelian, by Lemma \ref{l:lemm1} (ii), so $L^{(1)} \subseteq N_1(L)$. Repetition of this argument gives $N_i(L) = L^{(i)}$ for each $i \geq 0$.
\bigskip

If $F$ has characteristic zero, then every solvable Lie $A$-algebra over $F$ is metabelian, since $L^2$ is nilpotent. This is not the case, however, when $F$ is any field of characteristic $p > 0$, as the following example, which is taken from \cite[pages 52, 53]{jac}, shows.
\medskip

\begin{ex}\label{e:jac}
Let
$$ e = \left[ \begin{array}{rrrrrrr}
0 & 1 & 0 & . & . & . & 0\\
0 & 0 & 1 & 0 & . & . & 0\\
\vdots &  &  &  &  &  & \vdots \\
0 & . & . & . & . & 0 & 1 \\
1 & 0 & . & . & . & . & 0
\end{array} \right], \hspace{.2cm}
f = \left[ \begin{array}{rrrrr}
0 & 0 & 0 & \ldots & 0 \\
0 & 1 & 0 & \ldots & 0 \\
0 & 0 & 2 & \ldots & 0 \\
\vdots &  &  &  & \vdots \\
0 & 0 & 0 & \ldots & p-1
\end{array} \right],
$$ 
let $F$ be a field of prime characteristic $p$ and put $L = Fe + Ff + F^p$ with product $[a + {\bf x}, b + {\bf y}] = [a,b] + ({\bf x}b - {\bf y}a)$ for all $a, b \in Fe + Ff$, ${\bf x}, {\bf y} \in F^p$. 
\end{ex}
Then $L$ is a solvable Lie algebra and $L^2 = Fe + F^p$ is not nilpotent. Moreover, $F^p$ is a minimal ideal of $L$ so the maximal subalgebras are either isomorphic to $Fe + Ff$, which is solvable but not nilpotent, or of the form $F(\alpha e + \beta f) + F^p$ for some $\alpha, \beta \in F$ with $(\alpha, \beta) \neq (0,0)$. It is straightforward to calculate that the characteristic polynomial of $\alpha e + \beta f$ is $x^p - \beta^{p-1}x - \alpha^p$. This is never divisible by $x^2$ and is divisible by $x$ if and only if $\alpha = 0$. It follows that the nilpotent subalgebras of $L$ are one-dimensional, $Ff + F {\bf x}_1$ where ${\bf x}_1 = (1,0,0, \ldots, 0)$, or inside $F^p$; in particular, all of them are abelian so this is a Lie $A$-algebra. 
\par
Note that $L$ is also $\phi$-free but not elementary. For let $B = Fe + F^p$. Then it is easy to see that $F({\bf x}_1 + \dots + {\bf x}_p)$ (where ${\bf x}_i$ is the $i^{th}$ standard basis vector for $F^p$) is an ideal of $B$, and is, in fact, $\phi(B)$. Therefore this is an example of a Lie $A$-algebra that is not an $E$-algebra.
\bigskip

If $B$ is a subalgebra of $L$, the {\em centraliser} of $B$ in $L$ is $Z_L(B) = \{ x \in L : [x,B] = 0 \}$. We shall also need the following simple result. 
\bigskip

\begin{lemma}\label{l:nilrad}
Let $L$ be any solvable Lie algebra with nilradical $N$. Then $Z_L(N) \subseteq N$ 
\end{lemma}
\medskip
{\it Proof.} Suppose that $Z_L(N) \not \subseteq N$. Then there is a non-trivial abelian ideal $A/(N \cap Z_L(N)$ of $L/(N \cap Z_L(N)$ inside $Z_L(N)/(N \cap Z_L(N)$. But now $A^3 \subseteq [A, N] = 0$, so $A$ is a nilpotent ideal of $L$. It follows that $A \subseteq N \cap Z_L(N)$, a contradiction.

\section{Decomposition results}
\medskip
Here we have the basic structure theorems. First we see that $L$ splits over the terms in its derived series.
\bigskip

\begin{theor}\label{t:split}
Let $L$ be a solvable Lie $A$-algebra. Then $L$ splits over each term in its derived series. Moreover, the Cartan subalgebras of $L^{(i)}/L^{(i+2)}$ are precisely the subalgebras that are complementary to $L^{(i+1)}/L^{(i+2)}$ for $i \geq 0$.
\end{theor}
\medskip
{\it Proof.} Suppose that $L^{(n+1)} = 0$ but $L^{(n)} \neq 0$. First we show that $L$ splits over $L^{(n)}$. Clearly we can assume that $n \geq 1$. Let $C$ be a Cartan subalgebra of $L^{(n-1)}$ (see, for example, \cite[Corollary 4.4.1.2]{wint}) and let $L = L_0 \dot{+} L_1$ be the Fitting decomposition of $L$ relative to ad$C$. Then $L_1 = \cap_{k=1}^{\infty} L({\rm ad}C)^k \subseteq L^{(n)}$, and so $L_1$ is an abelian ideal of $L$. Also $L^{(n-1)} = L_1 \dot{+} L_0 \cap L^{(n-1)}$ and  $L_0 \cap L^{(n-1)} = (L^{(n-1)})_0 = C$, which is abelian. It follows that $L^{(n-1)}/L_1$ is abelian, whence $L^{(n)} \subseteq L_1$ and $L = L_0 \dot{+} L^{(n)}$.
\par
So we have that $L = L^{(n)} \dot{+} B$ where $B = L_0$ is a subalgebra of $L$. Clearly $B^{(n)} = 0$, so, by the above argument, $B$ splits over $B^{(n-1)}$, say $B = B^{(n-1)} \dot{+} D$. But then $L = L^{(n)} \dot{+} (B^{(n-1)} \dot{+} D) = L^{(n-1)} \dot{+} D$. Continuing in this way gives the desired result.
\bigskip

This gives us the following fundamental decomposition result.
\bigskip

\begin{coro}\label{c:decomp}
Let $L$ be a solvable Lie $A$-algebra of derived length $n+1$. Then
\begin{itemize}
\item[(i)] $L = A_{n} \dot{+} A_{n-1} \dot{+} \ldots \dot{+} A_0$ where $A_i$ is an abelian subalgebra of $L$ for each $0 \leq i \leq n$; and
\item[(ii)] $L^{(i)} = A_{n} \dot{+} A_{n-1} \dot{+} \ldots \dot{+} A_{i}$ for each $0 \leq i \leq n$
\end{itemize}
\end{coro}
\medskip
{\it Proof.} (i) By Theorem \ref{t:split} there is a subalgebra $B_n$ of $L$ such that $L = L^{(n)} \dot{+} B_n$. Put $A_n = L^{(n)}$. Similarly $B_n = A_{n-1} \dot{+} B_{n-1}$ where $A_{n-1} = (B_n)^{(n-1)}$. Continuing in this way we get the claimed result. Note, in particular, that it is apparent from the construction that $A_k \cap (A_{k-1} + ... + A_0) = 0$ for each $1 \leq k \leq n$, and that it is easy to see from this that the sum is a vector space direct sum.
\par
(ii) We have that $L^{(n)} = A_n$. Suppose that $L^{(k)} = A_n \dot{+} \ldots \dot{+} A_k$ for some $1 \leq k \leq n$. Then $L = L^{(k)} \dot{+} B_k$ and $A_{k-1} = B_k^{(k-1)}$ by the construction in (i). But now $L^{(k-1)} \subseteq L^{(k)} + B_k^{(k-1)} \subseteq L^{(k-1)}$, whence $L^{(k-1)} = A_n \dot{+} \ldots \dot{+} A_{k-1}$ and the result follows by induction.
\bigskip

Now we show that the result in Theorem \ref{t:a} (iii)(a) holds when $L$ is solvable without any restrictions on the underlying field. We say that $L$ is {\em monolithic} with {\em monolith} $W$ if $W$ is the unique minimal ideal of $L$. 
\bigskip

\begin{theor}\label{t:int}
Let $L$ be a solvable Lie $A$-algebra. Then $Z(L) \cap L^2 = 0$.
\end{theor}
\medskip
{\it Proof.} Let $L$ be a minimal counter-example and let $z \in Z(L) \cap L^2$. Put $Z(L) = U \dot{+} Fz$. Then $U$ is an ideal of $L$ and 
$$ U \neq z + U \in (Z(L) \cap L^2 + U)/U \subseteq Z(L/U) \cap (L/U)^2. $$
The minimality of $L$ implies that $U = 0$, so $Z(L) = Fz$. But now if $K$ is an ideal of $L$ which does not contain $Z(L)$, then $K \neq z + K \in Z(L/K) \cap (L/K)^2$ similarly, contradicting the minimality of $L$. It follows that $L$ is monolithic with monolith $Z(L)$.
\par
Now let $M$ be a maximal ideal of $L$. Then $Z(M) \cap M^2 = 0$ by the minimality of $L$, so $Z(L) \not \subseteq M^2$, whence $M^2 = 0$. It follows that $L = M \dot{+} Fx$ for some $x \in L$ and $M$ is abelian. Let $L = L_0 \dot{+} L_1$ be the Fitting decomposition of $L$ relative to ad$x$. Then $L_1 = \cap_{i=1}^{\infty} L({\rm ad}x)^i \subseteq M$, and $[L_0,L_1] \subseteq L_1$, so $L_1$ is an ideal of $L$. If $L_1 \neq 0$ then $Z(L) \subseteq L_1 \cap L_0 = 0$, a contradiction. Hence $L_1 = 0$ and ad\,$x$ is nilpotent. But then $L =M + Fx$ is nilpotent and hence abelian, and the result follows.
\bigskip 

Next we aim to show the relationship between ideals of $L$ and the decomposition given in Corollary \ref{c:decomp}. First we need the following lemma.
\bigskip

\begin{lemma}\label{l:ideal}
Let $L$ be a solvable Lie $A$-algebra of derived length $\leq n+1$, and suppose that $L = B \dot{+} C$ where $B = L^{(n)}$ and $C$ is a subalgebra of $L$. If $D$ is an ideal of $L$ then $D = (B \cap D) \dot{+} (C \cap D)$.
\end{lemma}
\medskip
{\it Proof.} Let $L$ be a counter-example for which dim$L$ + dim$D$ is minimal. Suppose first that $D^2 \neq 0$. Then $D^2 = (B \cap D^2) \dot{+} (C \cap D^2)$ by the minimality of $L$. Moreover, since
$$ L/D^2 = (B + D^2)/D^2 \hspace{.2cm} \dot{+} \hspace{.2cm} (C + D^2)/D^2 $$
we have
$$ D/D^2 = (B \cap D + D^2)/D^2 \hspace{.2cm} \dot{+} \hspace{.2cm} (C \cap D + D^2)/D^2 $$
whence 
$$D = B \cap D + C \cap D + D^2 = B \cap D \dot{+} C \cap D.$$
We therefore have that $D^2 = 0$. Similarly, by considering $L/B \cap D$, we have that $B \cap D = 0$. 
\par
Put $E = C^{(n-1)}$. Then $(D + B)/B$ and $(E + B)/B$ are abelian ideals of the Lie $A$-algebra $L/B$, and so 
$$ \left[ \frac{D + B}{B}, \frac{E + B}{B} \right] = \frac{B}{B}, $$
by Lemma \ref{l:lemm1} (ii), whence
$$ [D, E] \subseteq [D + B, E + B] \subseteq B \hbox{ and } [D,E] \subseteq B \cap D = 0;$$ that is, $D \subseteq Z_L(E)$. But $Z_L(E) = Z_B(E) + Z_C(E)$. For, suppose that $x = b + c \in Z_L(E)$, where $b \in B$, $c \in C$. Then $0 = [x,E] = [b,E] + [c,E]$, so $[b,E] = - [c,E] \in B \cap C = 0$. This implies that $Z_L(E) \subseteq Z_B(E) + Z_C(E)$. But the reverse inclusion is clear, so equality follows.
\par
Now $L^{(n-1)} \subseteq B + E \subseteq L^{(n-1)}$, so $B = L^{(n)} = (B + E)^2 = [B,E]$. Let $L^{(n-1)} = L_0 \dot{+} L_1$ be the Fitting decomposition of $L^{(n-1)}$ relative to ad\,$E$. Then $B \subseteq L_1$ so that $Z_B(E) \subseteq L_0 \cap L_1 = 0$, whence $D\subseteq Z_L(E) = Z_C(E) \subseteq C$ and the result follows.  
\bigskip

\begin{theor}\label{t:nz}
Let $L$ be a solvable Lie $A$-algebra of derived length $n+1$ with nilradical $N$, and let $K$ be an ideal of $L$ and $A$ a minimal ideal of $L$. Then, with the same notation as Corollary \ref{c:decomp}, 
\begin{itemize}
\item[(i)] $K = (K \cap A_n) \dot{+} (K \cap A_{n-1}) \dot{+} \ldots \dot{+} (K \cap A_0)$;
\item[(ii)] $N = A_n \oplus (N \cap A_{n-1}) \oplus \ldots \oplus (N \cap A_0)$;
\item[(iii)] $Z(L^{(i)}) = N \cap A_i$ for each $0 \leq i \leq n$; and
\item[(iv)] $A \subseteq N \cap A_i$ for some $0 \leq i \leq n$. 
\end{itemize} 
\end{theor}
\medskip
{\it Proof.} (i) We have that $L = A_n \dot{+} B_n$ where $A_n = L^{(n)}$ from the proof of Corollary \ref{c:decomp}. It follows from Lemma \ref{l:ideal} that $K = (K \cap A_n) + (K \cap B_n)$. But now $K \cap B_n$ is an ideal of $B_n$ and $B_n = A_{n-1} \dot{+} B_{n-1}$. Applying Lemma \ref{l:ideal} again gives $K \cap B_n = (K \cap A_{n-1}) \dot{+} (K \cap B_{n-1})$. Continuing in this way gives the required result.
\par
(ii) This is clear from (i), since $A_n = L^{(n)} = N \cap A_n$.
\par
(iii) We have that $L^{(i)} = L^{(i+1)} \dot{+} A_i$ from Corollary \ref{c:decomp}, and $Z(L^{(i)}) \cap L^{(i+1)} = 0$ from Theorem \ref{t:int}. Thus, using Lemma \ref{l:ideal},
$$ Z(L^{(i)}) = (Z(L^{(i)}) \cap L^{(i+1)}) + (Z(L^{(i)}) \cap A_i) = Z(L^{(i)}) \cap A_i \subseteq N \cap A_i. $$
It remains to show that $N \cap A_i \subseteq Z(L^{(i)})$; that is, $[N \cap A_i,L^{(i)}] = 0$. 
We use induction on the derived length of $L$. If $L$ has derived length one the result is clear. So suppose it holds for Lie algebras of derived length $\leq k$, and let $L$ have derived length $k+1$. Then $B = A_{k-1} + \dots + A_0$ is a solvable Lie $A$-algebra of derived length $k$, and, if $N$ is the nilradical of $L$, then $N \cap A_i$ is inside the nilradical of $B$ for each $0 \leq i \leq k-1$, so $[N \cap A_i, B^{(i)}] = 0$ for $0 \leq i \leq k-1$, by the inductive hypothesis. But $[N \cap A_i, A_k] = [N \cap A_i,L^{(k)}] \subseteq [N,N] = 0$, for $0 \leq i \leq k$, whence $[N \cap A_i,L^{(i)}] = [N \cap A_i,A_k + B^{(i)}] = 0$ for $0 \leq i \leq k$.
\par
(iv) We have $A \subseteq L^{(i)}$, $A \not \subseteq L^{(i+1)}$ for some $0 \leq i \leq n$. Now $[L^{(i)}, A] \subseteq [L^{(i)}, L^{(i)}] = L^{(i+1)}$, so $[L^{(i)}, A] \neq A$. It follows that $[L^{(i)}, A] = 0$, whence $A \subseteq Z(L^{(i)}) = N \cap A_i$, by (ii).
\bigskip

The final result in this section shows when two ideals of a Lie $A$-algebra centralise each other.
\bigskip

\begin{propo}\label{p:cent}
Let $L$ be a Lie $A$-algebra and let $B, D$ be ideals of $L$. Then $B \subseteq Z_L(D)$ if and only if $B \cap D \subseteq Z(B) \cap Z(D)$.
\end{propo}
\medskip
{\it Proof.} Suppose first that $B \subseteq Z_L(D)$. Then $[B \cap D,D] = 0 = [B \cap D,B]$, whence $B \cap D \subseteq Z(B) \cap Z(D)$.
\par
Conversely, suppose that $B \cap D \subseteq Z(B) \cap Z(D)$. Then $[B,D] \subseteq B \cap D \subseteq Z(B + D)$ which yields that $[B,D] \subseteq (B + D)^2 \cap Z(B + D) = 0$, by Theorem \ref{t:int}. Hence $B \subseteq Z_L(D)$.
\bigskip

\section{Strongly solvable Lie $A$-algebras}
\medskip
A Lie algebra $L$ is called {\it strongly solvable} if $L^2$ is nilpotent. Over a field of characteristic zero every solvable Lie algebra is strongly solvable. Clearly strongly solvable Lie $A$-algebras are metabelian so we would expect stronger results to hold for this class of algebras. First the decomposition theorem takes on a simpler form.
\bigskip

\begin{theor}\label{t:ss}
Let $L$ be a strongly solvable Lie $A$-algebra with nilradical $N$. Then $L = L^2 \dot{+} B$, where $L^2$ is abelian and $B$ is an abelian subalgebra of $L$, and $N = L^2 \oplus Z(L)$.
\end{theor}
\medskip
{\it Proof.} We have that $L = L^2 \dot{+} B$, where $B$ is an abelian subalgebra of $L$, by Theorem \ref{t:split}. Also, $L^2$ is nilpotent and so abelian. Moreover, $N = L^2 + N \cap B$ and $N \cap B = Z(L)$, by Theorem \ref{t:nz}.
\bigskip

Next we see that the minimal ideals are easy to locate.
\bigskip

\begin{theor}\label{t:min}
Let $L = L^2 \dot{+} B$ be a strongly solvable Lie $A$-algebra and let $A$ be a minimal ideal of $L$. Then 
\begin{itemize}
\item[(i)] $A \subseteq L^2$ or $A \subseteq B$;
\item[(ii)] $A \subseteq B$ if and only if $A \subseteq Z(L)$ (in which case dim $A = 1$); and
\item[(iii)] $A \subseteq L^2$ if and only if $[A,L] = A$.
\end{itemize} 
\end{theor}
\medskip
{\it Proof.} (i) and (ii) follow from Theorem \ref{t:nz} (iii) and (iv).
\par
(iii) Suppose that $A \subseteq L^2$. Then $[A,L] \neq 0$ from (ii), so $[A,L] = A$. The converse is clear.
\bigskip

\begin{coro}\label{c:phi}
Let $L$ be a strongly solvable Lie $A$-algebra. Then $L$ is $\phi$-free if and only if $L^2 \subseteq $ Asoc$L$.
\end{coro}
\medskip
{\it Proof.} Suppose first that $L$ is $\phi$-free. Then $L^2 \subseteq N = $ Asoc$L$, by \cite[Theorem 7.4]{frat}. 
\par
So suppose now that $L^2 \subseteq $ Asoc$L$. Then $L$ splits over Asoc$L$ by Theorem \ref{t:split}. But now $L$ is $\phi$-free by \cite[Theorem 7.3]{frat}.
\bigskip

Finally we can identify the maximal nilpotent subalgebras of $L$. First we need the following lemma.
\bigskip

\begin{lemma}\label{l:maxn}
Let $L$ be a metabelian Lie algebra, and let $U$ be a maximal nilpotent subalgebra of $L$. Then $U \cap L^2$ is an abelian ideal of $L$ and $L^2 = (U \cap L^2) \oplus K$ where $K$ is an ideal of $L$ and $[U,K] = K$.
\end{lemma}
\medskip
{\it Proof.} Let $L = L_0 \dot{+} L_1$ be the Fitting decomposition of $L$ relative to ad\,$U$. Then $L_1 = \cap_{i=1}^{\infty} L({\rm ad}\,U)^i \subseteq L^2$, and so $L_1$ is an abelian ideal of $L$. Moreover, $L^2 = (L_0 \cap L^2)\dot{+} L_1$ and 
$$[L,L_0 \cap L^2] = [L_0 + L_1,L_0 \cap L^2] \subseteq (L_0 \cap L^2) + L^{(2)} = L_0 \cap L^2,$$ 
so $L_0 \cap L^2$ is an ideal of $L$. It follows that $U + (L_0 \cap L^2)$ is nilpotent and so $L_0 \cap L^2 \subseteq U \cap L^2$. The reverse inclusion is clear. Finally put $K = L_1$. 
\bigskip 

\begin{theor}\label{t:maxn}
Let $L$ be a strongly solvable Lie $A$-algebra, and let $U$ be a maximal nilpotent subalgebra of $L$. Then $U = (U \cap L^2) \oplus (U \cap C)$ where $C$ is a Cartan subalgebra of $L$.
\end{theor}
\medskip
{\it Proof.} Put $U = (U \cap L^2) \oplus D$, so $D$ is an abelian subalgebra of $L$. Let $L = L_0 \dot{+} L_1$ be the Fitting decomposition of $L$ relative to ad\,$D$. As in Lemma \ref{l:maxn}, $L_1$ is an abelian ideal of $L$, so $L^2 = L_0^2 \oplus [L_0,L_1]$, whence $L_0 \cap L^2 = L_0^2$. 
\par
Now put $L^2 = (U \cap L^2) \oplus K$ as given by Lemma \ref{l:maxn}. Then 
$$K = [U,K] = [D,K] \hbox{ so } K \subseteq L_1 \hbox{ and } U \cap L^2 \subseteq L_0 \cap L^2.$$ 
Hence 
$$L_0^2 = L_0 \cap L^2 = (U \cap L^2) + (L_0 \cap L^2 \cap K) = U \cap L^2.$$
\par 
Next put $L_0 = L_0^2 \dot{+} E$ where $E$ is an abelian subalgebra of $L_0$. Then 
$$U = L_0 \cap U = L_0^2 \oplus (E \cap U) = (U \cap L^2) \oplus (E \cap U). \hspace{1cm} (*)$$
\par
Finally put $E = (E \cap L^2) \oplus C$ where $E \cap U \subseteq C$. Then 
$$L = L_1 + L_0 = L^2 + L_0 = L^2 + E = L^2 \dot{+} C $$ 
so $C$ is a Cartan subalgebra of $L$, by Theorem \ref{t:split}. Moreover, $E \cap U \subseteq C \cap U$, so (*) implies that 
$$C \cap U = (E \cap U) \oplus (C \cap U \cap L^2) = E \cap U,$$ 
since $C \cap L^2 = 0$. But now (*) becomes $U = (U \cap L^2) \oplus (U \cap C)$ where $C$ is a Cartan subalgebra of $L$, as claimed.   
\bigskip

\section{Monolithic solvable Lie $A$-algebras}
\medskip

Monolithic algebras play a part in the application of $A$-algebras to the study of residually finite varieties, so it seems worthwhile to investigate what extra properties they might have.
\bigskip

\begin{theor}\label{t:mon}
Let $L$ be a monolithic solvable Lie $A$-algebra of derived length $n+1$ with monolith $W$. Then, with the same notation as Corollary \ref{c:decomp}, 
\begin{itemize}
\item[(i)] $W$ is abelian;
\item[(ii)] $Z(L) = 0$ and $[L,W] = W$; 
\item[(iii)] $N = A_n = L^{(n)}$; 
\item[(iv)] $N = Z_L(W)$; and
\item[(v)] $L$ is $\phi$-free if and only if $W = N$.
\end{itemize} 
\end{theor}
\medskip
{\it Proof.} (i) Clearly $W \subseteq L^{(n)}$, which is abelian.
\par
(ii) If $Z(L) \neq 0$ then $W \subseteq Z(L) \cap L^2 = 0$, by Theorem \ref{t:nz}, a contradiction. Hence $Z(L) = 0$. It follows from this that $[L,W] \neq 0$, whence $[L,W] = W$.
\par
(iii) We have $N = A_n \oplus N \cap A_{n-1} \oplus \ldots \oplus N \cap A_0$ by Theorem \ref{t:nz}(i). Moreover, $N \cap A_i$ is an ideal of $L$ for each $0 \leq i \leq n-1$, by Theorem \ref{t:nz}(iii). But if $N \cap A_i \neq 0$ then $W \subseteq A_n \cap N \cap A_i = 0$ if $i \neq n$. This contradiction yields the result.
\par
(iv) We have that $L = N \dot{+} B$ for some subalgebra $B$ of $L$, by Theorem \ref{t:split} and (iii). Put $C = Z_L(W)$ and note that $N \subseteq C$. Suppose that $N \neq C$. Then $C = N \dot{+} B \cap C$. Choose $A$ to be a minimal ideal of $B \cap C$, so that $A$ is abelian, and let $L = L_0 \dot{+} L_1$ be the Fitting decomposition of $L$ relative to ad\,$A$. Then 
$$L_1 = \bigcap_{i=1}^{\infty} L ({\rm ad}A)^i \subseteq [[[L,A],A],A] \subseteq [[C,A],A] \subseteq [N+A,A] \subseteq N,$$ 
which is abelian. It follows that $L_1$ is an ideal of $L$ and so $L_1 =0$, since otherwise $W \subseteq L_1 \cap L_0 = 0$. This yields that $N + A$ is nilpotent and thus abelian, whence $A \subseteq Z_L(N) \subseteq N$, by Lemma \ref{l:nilrad}. This contradiction implies that $N = C$.
\par
(v) Clearly $W =$ Asoc$L$. Suppose first that $L$ is $\phi$-free. Then $W =$ Asoc$L = N$, by \cite[Theorem 7.4]{frat}. So suppose now that Asoc$L = W = N$. Then $L$ splits over Asoc$L$ by Theorem \ref{t:split} and (iii). But now $L$ is $\phi$-free by \cite[Theorem 7.3]{frat}.
\bigskip

Note that Example \ref{e:jac} is monolithic, so monolithic solvable $A$-algebras are not necessarily metabelian. However, when the Lie $A$-algebra is strongly solvable the situation is more straightforward.
\bigskip

\begin{theor}\label{t:ssmon}
Let $L$ be a monolithic strongly solvable Lie $A$-algebra. Then the maximal nilpotent subalgebras of $L$ are $L^2$ and the Cartan subalgebras of $L$ (that is, the subalgebras that are complementary to $L^2$.) 
\end{theor}
\medskip
{\it Proof.}
Let $U$ be a maximal nilpotent subalgebra of $L$ and let $W$ be the monolith of $L$. Then $L^2 = (U \cap L^2) \oplus K$ where $U \cap L^2, K$ are ideals of $L$ and $[U,K] = K$, by Lemma \ref{l:maxn}. Either $W \subseteq U \cap L^2$ and $K = 0$ or else $W \subseteq K$ and $U \cap L^2 = 0$.
\par
In the former case $N = L^2 \subseteq U$, by Theorem \ref{t:mon}. But then $U \subseteq Z_L(N) \subseteq N$, by Lemma \ref{l:nilrad}, so $U = L^2$. In the latter case $U$ is a Cartan subalgebra of $L$, by Theorem \ref{t:maxn}. 
\bigskip

Finally we give necessary and sufficient conditions for a monolithic algebra to be a strongly solvable Lie $A$-algebra. The next two results are essentially Lemma 3 of \cite{sh1}, though the proofs are somewhat different.
\bigskip

\begin{lemma}\label{l:aa}
Let $L = L^2 \dot{+} B$ be a metabelian Lie algebra, where $B$ is a subalgebra of $L$, and suppose that $[L^2,b] = L^2$ for all $b \in B$. Then $L$ is a strongly solvable $A$-algebra.
\end{lemma}
\medskip
{\it Proof.} Let $U$ be a maximal nilpotent subalgebra of $L$. We have $L^2 = (U \cap L^2) \oplus K$ where $K$ is an ideal of $L$ and $[U,K] = K$, by Lemma \ref{l:maxn}. Let $u = x + b \in U$, where $x \in L^2$, $b \in B$. Then $L^2 = [L^2,b] = [L^2,u]$, so $L^2 = L^2 ({\rm ad}\,u)^i$ for all $i \geq 1$. It follows that $L^2 = K$ from which $U^2 \subseteq U \cap L^2 = 0$ and $L$ is an $A$-algebra.
\bigskip

\begin{theor}\label{l:mona}
Let $L$ be a monolithic Lie algebra. Then $L$ is a strongly solvable $A$-algebra if and only if $L = L^2 \dot{+} B$ is metabelian, where $B$ is a subalgebra of $L$ and $[L^2,b] = L^2$ for all $b \in B$ (or, equivalently, ad\,$b$ acts invertibly on $L^2$).
\end{theor}
\medskip
{\it Proof.} Suppose first that $L$ is a strongly solvable $A$-algebra. Then $L = L^2 \dot{+} B$ is metabelian, where $B$ is a subalgebra of $L$, by Theorem \ref{t:split}. Let $b \in B$ and let $L = L_0 \dot{+} L_1$ be the Fitting decomposition of $L$ relative to ad\,$b$. It is easy to see, as in Lemma \ref{l:maxn}, that $L^2 = (L^2 \cap L_0) \dot{+} L_1$ and $L^2 \cap L_0$ and $L_1$ are ideals of $L$, so $L^2 = L^2 \cap L_0$ or $L^2 = L_1$ as $L$ is monolithic. The former implies that $[L^2,b] = 0$, but then $L^2$ and $Fb$ are ideals of $L$, which is impossible. It follows that $L^2 = L_1$, whence $[L^2,b] = L^2$. If $\theta = ({\rm ad}\,b)|_{L^2}$ then $L^2 = {\rm Ker}\,\theta \dot{+} {\rm Im}\,\theta$, so ${\rm Ker}\,\theta = \{0\}$ and $\theta$ is invertible.
\par
The converse follows from Lemma \ref{l:aa}.
\bigskip

\section{Solvable $A$-algebras over an algebraically closed field}
\medskip
First we need the following lemma.
\bigskip

\begin{lemma}\label{l:min} Let $L$ be a solvable Lie $A$-algebra over an algebraically closed field $F$ of characteristic $p > 0$. Let $K$ be an ideal of $L$, $A$ a minimal ideal of $L$ with $A \subseteq Z(K)$, and $N$ an ideal of $L$ containing $K$ and such that $N/K \subseteq N(L/K)$, the nilradical of $L/K$. Then dim $(N/Z_N(A)) \leq 1$.
\end{lemma}
\medskip
{\it Proof.} Put $\bar{L} = L/K$ and for each $x \in L$ write $\bar{x} = x + K$. Then $A$ is an irreducible $\bar{L}$-module, and hence an irreducible $U$-module, where $U$ is the universal enveloping algebra of $\bar{L}$. Let $\phi$ be the corresponding representation of $U$ and let $\bar{x} \in \bar{L}$, $n \in N$. Then $[[\bar{x},\bar{n}],\bar{n}] = \bar{0}$, whence $[\bar{x},\bar{n}^p] = 0$ and so $\bar{n}^p \in Z = Z(U)$. 
\par
Let $n_1, n_2 \in N$. Then $\bar{n}_1^p, \bar{n}_2^p \in Z$, so $\alpha_1 \bar{n}_1^p + \alpha_2 \bar{n}_2^p \in \hbox{ker}(\phi)$, for some $\alpha_1, \alpha_2 \in F$, since dim $\phi(Z) \leq 1$, by Schur's Lemma. Since $F$ is algebraically closed, there are $\beta_1, \beta_2 \in F$ such that $\alpha_1 = \beta_1^p, \alpha_2 = \beta_2^p$, so $(\beta_1 \bar{n}_1 + \beta_2 \bar{n}_2)^p = \beta_1^p \bar{n}_1^p + \beta_2^p \bar{n}_2^p \in \hbox{ker}(\phi)$, since $[\bar{n}_1,\bar{n}_2] = \bar{0}$. It follows that $A + F(\beta_1 n_1 + \beta_2 n_2)$ is a nilpotent subalgebra of $L$ and hence abelian. Thus $\beta_1 \bar{n}_1 + \beta_2 \bar{n}_2 \in \hbox{ker}(\phi)$ and so dim $\phi(\bar{N}) \leq 1$. Hence $Z_N(A)$ has codimension at most $1$ in $N$.
\bigskip

The following result was proved by Drensky in \cite{dren}. We include a proof since, as far as we know, no English translation of the proof has appeared.
\bigskip

\begin{theor}\label{t:dren} Let $L$ be a solvable Lie $A$-algebra over an algebraically closed field $F$. Then the derived length of $L$ is at most $3$.
\end{theor}
\medskip
{\it Proof.} First note that we can assume that the ground field is of characteristic $p > 0$, since otherwise $L$ is strongly solvable and so of derived length at most $2$. Suppose that $L$ has derived length $4$. 
\par
Let $A$ be a minimal ideal of $L$ contained in $L^{(3)}$. Then, putting $K = L^{(3)}$, $N = L^{(2)}$ in Lemma \ref{l:min}, we deduce that dim $(L^{(2)}/Z_{L^{(2)}}(A)) \leq 1$.
Suppose that dim $(L^{(2)}/Z_{L^{(2)}}(A)) = 1$. Put $S = L/Z_{L^{(2)}}(A)$. Then dim$(S^{(2)}) = 1$. It follows that $S/Z_L(S^{(2)}) \subseteq \hbox{Der}(S^{(2)})$ and so has dimension at most one, giving $[S^{(1)},S^{(2)}] = 0$. But now $S^{(1)}$ is nilpotent but not abelian. As $S$ must be an $A$-algebra, this is a contradiction. We therefore have that dim $(L^{(2)}/Z_{L^{(2)}}(A)) = 0$, whence $[A,L^{(2)}] = 0$.
\par
Now we can include $L^{(3)}$ in a chief series for $L$. So let $0 = A_0 \subset A_1 \subset \ldots \subset A_r = L^{(3)}$ be a chain of ideals of $L$ each maximal in the next. By the above we have $[A_i,L^{(2)}] \subseteq A_{i-1}$ for each $1 \leq i \leq r$. It follows that $L^{(2)}$ is a nilpotent subalgebra of $L$ and hence abelian. We infer that $L^{(3)} = 0$, a contradiction.
\par
Clearly if the derived length of $L$ is greater than $4$ then $L/L^{(4)}$ is a solvable Lie $A$-algebra of derived length $4$ and the same contradiction follows.
\bigskip

Using the above we can examine in more detail the structure of monolithic Lie $A$-algebras.
\bigskip

\begin{theor}\label{t:algmon} Let $L$ be a monolithic solvable Lie $A$-algebra of dimension greater than one over an algebraically closed field $F$, with monolith $W$. Then either
\begin{itemize}
\item[(i)] $L = L^2 \dot{+} Fb$ where $L^2$ is abelian and $L^2(\hbox{ad}\,b - \lambda 1)^k = 0$ for some $k > 0$ and some $0 \neq \lambda \in F$, and dim\,$W = 1$; or
\item[(ii)] $F$ has characteristic $p > 0$, dim\,$W = p$ and $L = L^{(2)} \dot{+} B$ where $L^{(2)}$ is abelian, $B = Fb + Fn$, $[n,b] = n$, $L^{(2)}(\hbox{ad}\,n - \lambda 1)^k = 0$ and $L^{(2)}((\hbox{ad}\,b)^p - ad\,b - \mu^p 1)^k = 0$ for some $k > 0$ and some $0 \neq \lambda, \mu \in F$. 
\end{itemize}
\end{theor}
\medskip
{\it Proof.} Suppose first that $L$ is strongly solvable. Then $L = L^2 \dot{+} B$ where $L^2$ is abelian, $B$ is an abelian subalgebra and $W \subseteq L^2$. Now $W$ is an irreducible $B$-module and so one dimensional, by \cite[Section 1.5, Lemma 5.6]{sf}. Now $L/Z_L(W)$ is isomorphic to a subalgebra of Der$(W)$ and so $N = Z_L(W)$ has codimension at most one in $L$. It follows that $L$ is abelian (and hence one dimensional) or dim\,$B = 1$ and $N = L^2$. Decompose $L^2$ into ad\,$B$-invariant subspaces. Each is an ideal of $L$ and so there can be only one. It follows that $L^2(\hbox{ad}\,b - \lambda 1)^k = 0$ for some $k > 0$ and some $0 \neq \lambda \in F$, where $B = Fb$, giving case (i).
\par
So suppose now that $L^2$ is not nilpotent. Then $F$ has characteristic $p > 0$, $L$ has derived length $3$ and $W \subseteq L^{(2)}$. Let $N/L^{(2)}$ be the nilradical of $L/L^{(2)}$. Then applying Lemma \ref{l:min} with $K = L^{(2)}$ we see that dim\,$(N/Z_N(W)) \leq 1$. But $Z_L(W) = L^{(2)} \subseteq N$ by Theorem \ref{t:mon}(iv), so $Z_N(W) = L^{(2)}$. As $L^{(1)} \subseteq N$ we cannot have $N = L^{(2)}$, so dim\,$(N/L^{(2)}) = 1$; say $N = L^{(2)} \dot{+} Fn = L^{(1)}$.
\par
Let $L = L_0 \dot{+} L_1$ be the Fitting decomposition of $L$ relative to ad\,$n$. Then $L^{(2)} = [L^{(1)}, L^{(1)}] = [L^{(2)},n]$, so $L^{(2)} \subseteq L_1 \subseteq L^{(2)}$. Put $B = L_0$, so $B$ is a subalgebra of $L$ containing $n$ such that $L = L^{(2)} \dot{+} B$, and let $C = Z_B(Fn)$. Now $Fn = B \cap L^{(1)}$, giving that $Fn$, and hence $C$, is an ideal of $B$. Moreover, as $C^2 = Fn$, $C$ is a nilpotent ideal of $B$, and so $C = Fn$. It follows that $B/Fn = B/C$ has dimension at most one, and so dim\,$B \leq 2$. As $B$ is not abelian we have $B = Fn + Fb$ where $[n,b] = n$. This algebra has a unique $p$-map making it into a restricted Lie algebra: namely $b^{[p]} = b, n^{[p]} = 0$ (see \cite{sf}). We can decompose $L^{(2)} = \oplus_{\lambda,S} V_{\lambda,S}$ where $\lambda \in (Fn)^{*}, S \in B^{*}$ and 
\[
V_{\lambda,S} = \{x \in L^{(2)} : x(\hbox{ad}\,n - \lambda(n) 1)^k = 0 \hbox{ and } x((\hbox{ad}\,b)^p - ad\,b - S(b)^p 1)^k = 0 \}
\] 
by \cite[page 236]{sf}. As each $V_{\lambda,S}$ is an ideal of $L$ there can be only one of them. The fact that dim\,$W = p$ follows from \cite[Example 1, page 244]{sf}, so we have case (ii). 
\bigskip

\begin{coro} \label{c:algmon} If, in addition to the hypotheses of Theorem \ref{t:algmon}, $L$ is also $\phi$-free, then either
\begin{itemize}
\item[(i)] $L$ is two-dimensional non-abelian; or
\item[(ii)] $F$ has characteristic $p > 0$ and $L$ is isomorphic to the algebra in Example \ref{e:jac}. 
\end{itemize}
\end{coro}
\medskip
{\it Proof.} Case (i) follows from Theorem \ref{t:algmon} (i) because $W = L^2$ by Theorem \ref{t:mon}. If case (ii) of Theorem \ref{t:algmon} holds, then $F$ has characteristic $p > 0$, dim\,$W = p$ and $L = W \dot{+} B$ where $W$ is abelian, $B = Fb + Fn$ and $[n,b] = n$, using Theorem \ref{t:mon}. Let $\lambda$ be an eigenvalue for $(\hbox{ad}\,b)|_W$, so $[w,b] = \lambda\,w$ for some $w \in W$. Then $[w(\hbox{ad}\,n)^i,b] = (\lambda + i) w (\hbox{ad}\,n)^i$ for every $i$, so putting $w_i = w (\hbox{ad}\,n)^i$ we see that $Fw_0 + \dots + Fw_{p-1}$ is $B$-stable and hence equal to $W$. We then have $[w_i,b] = (\lambda + i) w_i$, $[w_i,n] = w_{i+1}$ (indices modulo $p$). But now the characteristic polynomial of ad\,$(b + \alpha\, n)$ is $(x - \lambda)^p - (x - \lambda) - \alpha^p$ and this is divisible by $x$ precisely when $\alpha^p = \lambda - \lambda^p$. It follows that by choosing $\alpha$ satisfying this equation and replacing $b$ by $b + \alpha\, n$ we can take $\lambda = 0$. This gives the algebra in Example \ref{e:jac}.   
\bigskip

Note: alternatively, it can be deduced that $W$ has the form claimed in (ii) by using \cite[Example 1, page 244]{sf}.
\par
Finally we seek describe the structure of $\phi$-free solvable Lie $A$-algebras over an algebraically closed field. The strongly solvable ones are easily described.
\bigskip

\begin{theor}\label{t:ssphifree} Let $L$ be a $\phi$-free strongly solvable Lie $A$-algebra over an algebraically closed field $F$. Then 
$$L = \sum_{i=1}^{m}Fa_i + \sum_{i=1}^{n}Fb_i \,\,\, {\rm where} \,\,\, [a_i,b_j]=\lambda_{ij}a_i $$
for all $1\leq i\leq m, 1\leq j\leq n$, other products being zero.
\end{theor}
\medskip
{\it Proof.} If $L$ is strongly solvable then it is elementary, by \cite[Theorem 2.5]{tv1}, and hence as described, by \cite[Theorem 3.2 (2)]{tv1}. (The restriction on the characteristic in that result is not required for the solvable case.)
\bigskip

The $\phi$-free solvable Lie $A$-algebras that are not strongly solvable are more complicated.
\bigskip

\begin{theor}\label{t:sphifree} Let $L$ be a $\phi$-free solvable Lie algebra, over an algebraically closed field $F$, that is not strongly solvable. Then $L$ is an $A$-algebra if and only if the following conditions are satisfied:
\begin{itemize}
\item[(i)] $L = L^{(2)} \dot{+} C \dot{+} B$, where $B$, $C$ are abelian subalgebras of $L$ and $L^{(1)} = L^{(2)} \dot{+} C$;
\item[(ii)] $B \dot{+} C$ is a strongly solvable $\phi$-free Lie $A$-algebra (and hence given by Theorem \ref{t:ssphifree});
\item[(iii)] $L^{(2)} = A_1 \oplus \ldots \oplus A_n$, where $A_i$ is a minimal ideal of $L$ of dimension $p$ for each $1 \leq i \leq n$; and
\item[(iv)] for each $1 \leq i \leq n$, there exists $c_i \in C$, $b_i \in B$ and a basis $a_{i1}, \ldots, a_{ip}$ for $A_i$ such that $C = Z_C(A_i) \oplus Fc_i$, $B = Z_B(A_i) \oplus Fb_i$, $[c_i, b_i] = c_i$, $[a_{ij}, c_i] = a_{i(j+1)}$ (indices modulo $p$) and $[a_{ij}, b] = (\lambda_i + j) a_{ij}$ for $1 \leq j \leq p$ and some $\lambda_i \in F$.
\end{itemize}  
\end{theor}
\medskip
{\it Proof.} Suppose first that $L$ is a $\phi$-free solvable Lie $A$-algebra that is not strongly solvable. Then $F$ has characteristic $p > 0$, $L = L^{(2)} \dot{+} C \dot{+} B$ where $L^{(2)}$ is abelian, $B$, $C$ are abelian subalgebras of $L$ and $L^{(1)} = L^{(2)} \dot{+} C$, by Theorem \ref{t:dren} and Corollary \ref{c:decomp}; this is (i). Moreover, $L^{(2)} \subseteq N(L) =$ Asoc $L$ , by \cite[Theorem 7.4]{frat}, so we can put $L^{(2)} = A_1 \oplus \ldots \oplus A_n$, where $A_i$ is a minimal ideal of $L$ for each $1 \leq i \leq n$. Put $L_i = A_i \dot{+} C \dot{+} B$. Then $L_i^{(1)} = A_i \dot{+} C$ and $L_i^{(2)} = A_i$, so $[C, B] = C$ and $[A_i, C] = A_i$. 
\par
Suppose first that dim\,$A_i = 1$. Then dim\,$L_i/Z_{L_i}(A_i) \leq 1$. If $A_i \dot{+} C = L_i^{(1)} \subseteq Z_{L_i}(A_i)$ then $A_i = [A_i, C] = 0$, a contradiction; so $C \not \subseteq Z_{L_i}(A_i)$. But $Z_{L_i}(A_i) = A_i \dot{+} (Z_{L_i}(A_i) \cap C) \dot{+} (Z_{L_i}(A_i) \cap B)$, by Theorem \ref{t:nz}, so $Z_{L_i}(A_i) \cap B = B$, giving $[A_i, B] = 0$. Hence $A_i = [A_i, C] = [A_i, [C, B]] \subseteq [C, [B, A_i]] + [B, [A_i, C]] = 0$, a contradiction again. It follows that dim\,$A_i \neq 1$. 
\par
Put $Z = Z_C(A_i) \dot{+} Z_B(A_i)$ and $\bar{L}_i = L_i/Z$. We claim that $\bar{L}_i$ is monolithic and $\phi$-free.
\par
Let $\bar{D} = D/Z$ be an ideal of $\bar{L}_i$ and suppose that $\bar{A}_i = (A_i + Z)/Z \not \subseteq \bar{D}$. Then $[A_i, D] \subseteq A_i \cap D = 0$, so $D \subseteq Z_{L_i}(A_i) \cap D = (A_i + Z) \cap D = (A_i \cap D) + Z = Z$. It follows that $\bar{L}_i$ is monolithic with monolith $\bar{A}_i$. Let $\bar{U} = U/Z$ be the nilradical of $\bar{L}_i$. Then $\bar{A}_i \subseteq \bar{U}$, so $A_i \subseteq U$ and $[A_i, U] \subseteq A_i \cap Z = 0$. This yields that $U \subseteq A_i + Z$, whence $\bar{U} = \bar{A}_i$. Theorem \ref{t:mon}(v) now implies that $\bar{L}_i$ is $\phi$-free. 
\par
Next put $D = C \dot{+} B$. Then $D/Z$ is two dimensional, by Corollary \ref{c:algmon}, and so $\phi$-free, whence $\phi(D) \subseteq Z \cap C = Z_C(A_i)$ for each $1 \leq i \leq n$. It follows that $\phi(D)$ is an ideal of $L$ and hence that $\phi(D) \subseteq \phi(L) = 0$, by \cite[Lemma 4.1]{frat}. This establishes (ii).
\par
Now $D$ is elementary, by \cite[Theorem 2.5]{tv1}, and so splits over each of its ideals, by Lemma 2.3 of \cite{elem}. This yields that $D = Z \dot{+} E$ for some subalgebra $E$ of $D$, whence $A_i \dot{+} E \cong \bar{L}_i$ has the form given in Corollary \ref{c:algmon}. Assertions (iii) and (iv) now follow.
\par
Now suppose that conditions (i)-(iv) are satisfied. Adopting the same notation as above we have that $L_i/Z$ is an $A$-algebra, by (iv), and that $L_i/A_i$ is an $A$-algebra, by (ii). It follows that $L_i$ is an $A$-algebra, by Lemma \ref{l:lemm2}. As this is true for each $1 \leq i \leq n$ repeated use of Lemma \ref{l:lemm2} yields that $L$ is an $A$-algebra.

\end{document}